# Multiply Periodic Splines

Zhao Guohui

School of Mathematical Sciences, Dalian University of Technology

Dalian City, Liaoning Province, P.R. China, 116024

Zhao Guohui, School of Mathematical Sciences, Dalian University of Technology,

No.2 Linggong Road, Ganjingzi District, Dalian City, Liaoning Province,

P.R. China, 116024, ghzhao@dlut.edu.cn, 86-0411-84708351-8315

# Multiply Periodic Splines

Zhao Guohui

School of Mathematical Sciences, Dalian University of Technology,

No. 2 Linggong Road, Ganjingzi District, 116024 Dalian Liaoning, P.R. China

**Abstract**

Spline functions have long been used in numerical solution of differential equations. Recently it revives as isogeometric analysis, which offers integration of finite element analysis and NURBS based CAD into a single unified process. Usually many NURBS pieces are needed to build geometrically continuous CAD models. In this paper, we introduce some multiply periodic splines defined on hyperbolic disc. A single piece of such splines is enough to build complex CAD models. Multiresolution analysis on surfaces of high genus built from such splines can be carried out naturally. CAD and FEA are integrated directly on such models. It is difficult to derive such splines, only a theoretical framework is presented, together with some simple examples. Rigorous derivation and construction of B-splines will be given in future papers.

## 1. Introduction

Isogeometric analysis integrates finite element analysis and tensor-product splines based, especially NURBS based CAD into a unified process, where triangulation for FEA is not needed. Usually a CAD model consists of many NURBS pieces, smoothly connected but independently defined. Most of the time, the parametric domains of all the pieces can not unified into a single domain, which results in difficulty for isogeometric analysis. In this paper, multiply periodic splines are introduced. Usually a single piece of such splines is enough to build complex CAD models, which facilitates isogeometric analysis.

Multiply periodic splines are defined on Klein's disk. The method to construct tensor-product splines is not enough. In section two, we introduce general bivariate splines, the method for studying general bivariate splines can be generalized to study multiply periodic splines. In section three, Klein disk is covered. In section four, multiply periodic partition and multiply periodic splines are defined, and multiply periodic smooth connection conditions are given. Some simple examples are given in section five, and summary comes at last.

## 2. Multivariate splines

Multivariate spline theory has been studied extensively. It is used in approximation, CAGD, FEM, etc. For the sake of simplicity, we restrict to the bivariate case. Given a polygonal domain $\Omega$, we partition it with irreducible algebraic curves into finite cells.

The partition is denoted by $\Delta$. A function $f(x)$ defined on $\Omega$ is a spline of degree $n$ and smooth order $r$ if $f \in c^r(\Omega)$, and on each cell $f(x)$ is a polynomial of degree $n$. The set of all splines of degree $n$ and smooth order $r$ is a linear space and is denoted by $s_n^r(\Omega, \Delta)$.

Several theorems of multivariate splines are needed in this paper and listed below, The interested readers may refer to [4] for their proofs.

**Theorem 1** A piecewise polynomial $p(x,y)$ defined on $\Omega$ with respect to $\Delta$ is a spline of degree $n$ and smooth order $r$ if and only if for any adjacent cells $\Delta_i$ and $\Delta_j$ with common partition curve $l_{ij}$

$$p_{|\Delta_i} - p_{|\Delta_j} = l_{ij}^{r+1} q_{ij}$$

where $q_{ij}$ is a polynomial and is called the smooth cofactor on $l_{ij}$.

**Theorem 2** Around an interior vertex of $\Delta$, smooth cofactors on partitioning curves concurrent at this vertex satisfy the following equation: $\sum l_{ij}^{r+1} q_{ij} \equiv 0$. This is called the conformality equation around this vertex.

**Theorem 3** The solution space of the conformality equation $\sum_{l=1}^{N} (\alpha_l x + \beta_l y)^{r+1} q_l \equiv 0, \alpha_i \beta_j \neq \alpha_j \beta_i, i \neq j$ has dimension

$$\frac{1}{2}(n - r - [(r+1)/(N-1)])_+ \times ((N-1)n - (N+1)r + (N-3) + (N-1)[(r+1)/(N-1)])$$

where $[x]$ denotes the largest integer $\leq x$.

## 3. Klein disk

Klein disk is a model of hyperbolic geometry in which points are represented by the points in the interior of the unit disk and lines are represented by the chords. Poincaré disk is another model of hyperbolic geometry. An advantage of the Poincaré disk is that it is conformal (circles and angles are not distorted); a disadvantage is that lines of the geometry are circular arcs orthogonal to the boundary circle of the disk.

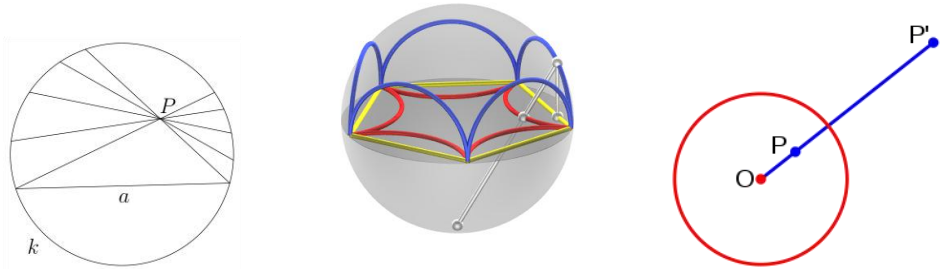

If P is a point a distance $u$ from the centre of the unit circle in the Poincaré disk, then the corresponding point of the Klein disk is a distance of $s$ on the same radius:

$$s = \frac{2u}{1+u^2}$$

Given two points $u$ and $v$ in the Poincaré disk which do not lie on a diameter, we can solve for the circle orthogonal to the boundary circle passing through both points, and obtain

$$x^2 + y^2 + \frac{u_2(v_1^2+v_2^2) - v_2(u_1^2+u_2^2) + u_2 - v_2}{u_1 v_2 - u_2 v_1} x$$
$$+ \frac{v_1(u_1^2+u_2^2) - u_1(v_1^2+v_2^2) + v_1 - u_1}{u_1 v_2 - u_2 v_1} y + 1 = 0$$

In the plane, the inverse of a point $P$ with respect to a reference circle with center $O$ and radius $r$ is a point $P'$ lying on the ray from $O$ through $P$ such that

$$OP \times OP' = r^2$$

Isometric transformations of Poincaré disk form a group, which is generated by reflections about geodesics in Poincaré disk. Such a reflection is represented in the model as a circle inversion about the circle that represents the geodesic. Rotations and translations can be represented as a combination of two reflections about different geodesics. In the case of rotations, the two intersect, while in the case of translations, they do not. The orientation preserving isometric transformations consist of Möbius transformations of the form

$$T(z) = \lambda \frac{z-a}{\bar{a}z - 1}, |\lambda| = 1, |a| < 1$$

Similarly, a geodesic in Klein disk is of the form $ux + vy = r, u^2 + v^2 = 1, 0 < r < 1$.

Isometric transformations of Klein disk form a group, which is generated by reflections about chords in Klein disk. Such a reflection is represented as a projective mapping. Rotations and translations can be represented as a combination of two reflections about different chords. In the case of rotations, the two intersect, while in the case of translations, they do not.

The projective plane $PR^2$ is defined as the quotient space of $R^3 \setminus \{(0,0,0)\}$ under the equivalent relation $(x, y, z) \sim (\lambda x, \lambda y, \lambda z), \lambda \neq 0$. Any $(x, y, z) \neq (0,0,0)$ is called the homogeneous coordinates of the projective plane.

For any nonsingular $3 \times 3$ matrix $A$, the linear mapping $x \to Ax$ induces a projective mapping $x \in PR^2 \to Ax \in PR^2$ by means of the homogeneous coordinates, which corresponds to a fractional linear mapping of the Euclidean plane:

$$x' = \frac{a_{11}x + a_{12}y + a_{13}}{a_{31}x + a_{32}y + a_{33}}, y' = \frac{a_{21}x + a_{22}y + a_{23}}{a_{31}x + a_{32}y + a_{33}}$$

The projective mappings which leave the unit circle invariant are called automorphic collineations, which form a group under mapping composition. Given a Klein line $l$ and a point $P$ in the Klein disk, the reflection point $P'$ of $P$ across $l$ can be obtained according to the property that $l$ is the perpendicular bisector of $\overrightarrow{PP'}$. For the special case that $l$ is defined by $x = \cos\theta$, the projective matrix $A$ can be obtained

$$A = \begin{pmatrix} 1+\cos^2\theta & 0 & -2\cos\theta \\ 0 & -\sin^2\theta & 0 \\ 2\cos\theta & 0 & -1-\cos^2\theta \end{pmatrix}$$

**4. Multiply periodic splines**

According to the uniformization theorem, all compact Riemann surfaces except sphere and torus are quotients of the Poincaré disk under some discrete transformation groups. It is called a Fuchsian group. We want to design and analyze high genus surfaces, we may construct multiply periodic splines on hyperbolic disk to build high genus surfaces. For a Fuchsian group of the unit disk, under which the quotient is a high genus surface, we select a partition of a fixed fundamental domain of the Fuchsian group so that when we transfer the partition to the whole disk using the Fuchsian group, we get a partition of the disk. Such a partition is called a multiply periodic partition. If cells of the partition are all triangles, it is called a triangulation. A function on the unit disk is called a multiply periodic function under the Fuchsian group if it is stable under the action of the Fuchsian group. A multiply periodic function can be regarded as a function defined on the quotient surface. A function defined on a fundamental domain can be extended to a multiply periodic function on the unit disk by the Fuchsian group action.

The geodesic of the Poincaré disk is a circle orthogonal to the unit circle, the boundary of a fundamental domain usually consists of consecutive arcs; while the geodesic of the Klein disk is a chord, a fundamental domain can be selected so that its boundary is a convex Euclidean polygon. So, we decide to study multiply periodic splines on Klein disk.

Given a multiply periodic partition of Klein disk, a $c^r$ function defined on Klein disk is an $s_n^r$ multiply periodic spline if it is both a $c^r$ multiply periodic function on the disk and an $s_n^r$ spline on a fundamental domain. To construct $s_n^r$ multiply periodic splines, we first construct an $s_n^r$ spline on the fundamental domain, then extend it to the whole disk to obtain a $c^r$ multiply periodic function. The function obtained is an $s_n^r$ multiply periodic spline if and only if it is continuous up to order $r$ across the boundary of the fundamental domain. This requirement results in conformality conditions similar to those of multivariate splines.

Suppose that the Fuchsian group $\Gamma$ is generated by $g_i, i=1,\cdots,k$. Two polynomials $p_i$ and $p_j$ are defined on two cells $\Delta_i$ and $\Delta_j$ respectively on the fundamental domain of the multiply periodic partition. Cell $\Delta_i$ is transferred to $g\Delta_i$ by a generator of the Fuchsian group. Suppose that $g\Delta_i$ and $\Delta_j$ are adjacent. Their common edge is $l_{ij}$. The function on $g\Delta_i$ is $gp_i(x) = p_i(gx)$, which is a rational function, written as $gp_i(x) = \frac{u(x)}{v(x)}$ with $v$ being a known polynomial determined by the generator $g$.
The two functions $p_j$ and $gp_i$ are continuous up to order $r$ across $l_{ij}$ if and only if $p_j v - u = l_{ij}^{r+1} q$ where $q$ is a polynomial, which is a system of linear equations in the coefficients of $p_j, u$ and $q$.

To make the idea clear, I present an example below. As the Poincaré disk is often used to introduce hyperbolic geometry, I take Bolza surface as the example.

Let us now recall the definition of the Bolza surface. Consider the regular hyperbolic octagon, with angles equal to $\pi/4$ and corners at $2^{-1/4} e^{i(\frac{\pi}{8}+\frac{k\pi}{4})}, k=0,\cdots,7$. The four hyperbolic translations $a$, $b$, $c$, and $d$ that identify opposite sides of the octagon generate a Fuchsian group. The Bolza surface is defined as the quotient of the Poincaré disk under the action of the Fuchsian group. Suppose $g_0, g_1, \cdots, g_7 = a, \bar{b}, c, \bar{d}, \bar{a}, b, \bar{c}, d$, where $\bar{g}$ denotes the inverse of $g$, the matrices of $g_k, k = 0, \cdots, 7$ are

$$g_k = \begin{pmatrix} \beta^2 & e^{ik\pi/4}\sqrt{2}\beta \\ e^{-ik\pi/4}\sqrt{2}\beta & \beta^2 \end{pmatrix}, \text{ where } \beta = \sqrt{1+\sqrt{2}}.$$

These translations are compositions of two 180°-rotations around the midpoint of a boundary edge and around the center of the octagon.

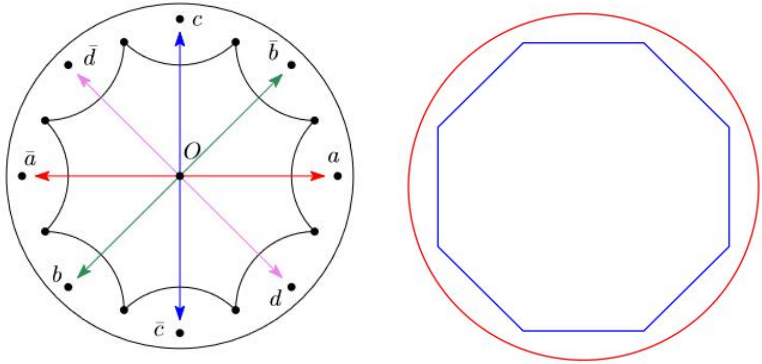

According to the relation between the Poincaré disk and the Klein disk, the corresponding regular octagon on the Klein disk has corners at

$$(2\sqrt{2}-2)2^{1/4}e^{i(\frac{\pi}{8}+\frac{k\pi}{4})}, k = 0, \cdots, 7$$

The automorphism of the Poincaré disk is of the form $T(z) = \dfrac{az+b}{\bar{b}z+\bar{a}}, a\bar{a} - b\bar{b} = 1$, where $a = a_1 + ia_2, b = b_1 + ib_2$. Since

$$\frac{2T(z)}{1+|T(z)|^2} = \frac{2\dfrac{az+b}{\bar{b}z+\bar{a}}}{1+\dfrac{az+b}{\bar{b}z+\bar{a}}\dfrac{\overline{az}+\bar{b}}{b\bar{z}+a}} = 2\frac{ab+\dfrac{aaz+bb\bar{z}}{1+|z|^2}}{|a|^2+|b|^2+2\dfrac{a\bar{b}z+\bar{a}b\bar{z}}{1+|z|^2}},$$

the corresponding automorphism of the Klein disk is of the form

$$\begin{pmatrix} a_1^2 - a_2^2 + b_1^2 - b_2^2 & 2b_1b_2 - 2a_1a_2 & 2a_1b_1 - 2a_2b_2 \\ 2a_1a_2 + 2b_1b_2 & a_1^2 - a_2^2 - b_1^2 + b_2^2 & 2a_1b_2 + 2a_2b_1 \\ 2a_1b_1 + 2a_2b_2 & 2a_1b_2 - 2a_2b_1 & a_1^2 + a_2^2 + b_1^2 + b_2^2 \end{pmatrix}$$

From the above formula, the generators $g_0, g_1, \cdots, g_7 = a, \bar{b}, c, \bar{d}, \bar{a}, b, \bar{c}, d$ of the corresponding Fuchsian group of the regular Klein octagon have the following form:

$$\begin{pmatrix} 3+2\sqrt{2}+(2+2\sqrt{2})\cos\dfrac{k\pi}{2} & (2+2\sqrt{2})\sin\dfrac{k\pi}{2} & (2+2\sqrt{2})^{3/2}\cos\dfrac{k\pi}{4} \\ (2+2\sqrt{2})\sin\dfrac{k\pi}{2} & 3+2\sqrt{2}-(2+2\sqrt{2})\cos\dfrac{k\pi}{2} & (2+2\sqrt{2})^{3/2}\sin\dfrac{k\pi}{4} \\ (2+2\sqrt{2})^{3/2}\cos\dfrac{k\pi}{4} & (2+2\sqrt{2})^{3/2}\sin\dfrac{k\pi}{4} & 5+4\sqrt{2} \end{pmatrix}$$

where $k = 0, \cdots, 7$.

## 5. $s_1^0$ multiply periodic splines

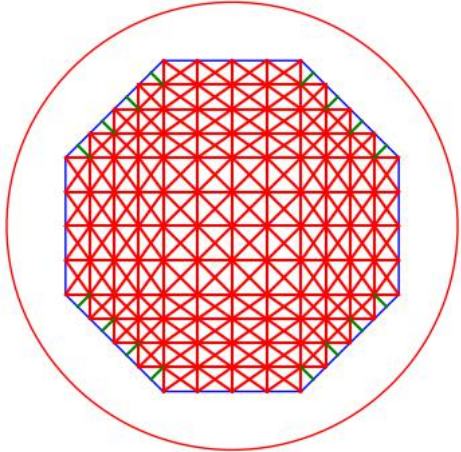

The regular octagon on the Klein disk corresponding to Bolza surface is triangulated in the above figure. It is easy to construct $s_1^0$ multiply periodic splines on this triangulation. These splines can be used to construct surfaces of genus two, and isogeometric analysis can be carried out on such surfaces.

## 6. Summary

Multiply periodic splines are introduced. B-splines and dimension of such spline spaces will be discussed later.